\documentclass[journal]{IEEEtran}

\usepackage{amsmath}
\usepackage{amsbsy}
\usepackage{graphicx}
\usepackage{subfigure}
\usepackage{cite}
\usepackage{algorithm}
\usepackage{algorithmic}
\usepackage{amsfonts}
\usepackage{amssymb}
\usepackage{graphicx}

\newtheorem{theorem}{Theorem}

\newtheorem{assumption}{Assumption}

\ifodd 1

\newcommand{\com}[1]{{#1}}
\newcommand{\res}[1]{\textbf{\color{magenta} (RESPONSE: #1)}}
\newcommand{\del}[1]{}
\else

\newcommand{\com}[1]{}
\newcommand{\res}[1]{}
\newcommand{\del}[1]{}
\fi

\begin{document}
\title{Enhanced Secondary Frequency Control \\via Distributed Peer-to-Peer Communication}

\author{Chenye~Wu,~\IEEEmembership{Member,~IEEE,}
        Soummya Kar,~\IEEEmembership{Member,~IEEE,}
        and~Gabriela Hug,~\IEEEmembership{Member,~IEEE} \vspace{-0.4cm}
\thanks{C.\,Wu, S. Kar and G. Hug are with the Department of Electrical and Computer Engineering, Carnegie Mellon University, Pittsburgh, PA 15213, USA, e-mails: chenyewu@andrew.cmu.edu, soummyak@andrew.cmu.edu, and ghug@ece.cmu.edu. This research was supported in part by Pennsylvania Infrastructure Technology Alliance and Carnegie Mellon University Scott Institute.} }

\maketitle

\begin{abstract}
Distributed generation resources have become significantly more prevalent in the electric power system over the past few years. This warrants reconsideration on how the coordination of generation resources is achieved. In this paper, we particularly focus on secondary frequency control and how to enhance it by exploiting peer-to-peer communication among the resources. We design a control framework based on a consensus-plus-global-innovation approach, which guarantees bringing the frequency back to its nominal value. The control signals of the distributed resources are updated in response to a global innovation corresponding to the ACE signal, and additional information exchanged via communication among neighboring resources. We show that such a distributed control scheme can be very well approximated by a PI controller and can stabilize the system. Moreover, since our control scheme takes advantage of both the ACE signal and peer-to-peer communication, simulation results demonstrate that our control scheme can stabilize the system significantly faster than the AGC framework. Also, an important feature of our scheme is that it performs $c\epsilon$-close to the centralized optimal economic dispatch, where $c$ is a positive constant depending only on the cost parameters and the communication topology and $\epsilon$ denotes the maximum rate of change of overall system.
\end{abstract}



\vspace{-0.2cm}
\section{Introduction}
\label{intro}

The electric power system is facing new challenges but also opportunities as more and more distributed resources such as  distributed generation units, storage systems, and demand response participants are integrated into the system. The discrepancy between this increasingly distributed infrastructure and the traditional centralized control structure may lead to unnecessary inefficiencies and reduction in reliability in the electric energy supply. On the other hand, the trend of incorporating more communication technology to the system allows for an enhancement of the control structure by leveraging the capability of direct communication and coordination between distributed resources.

\vspace{-0.2cm}
\subsection{Scope of the Work}

Given the opportunities brought by communication technology, we aspire to study how communications between distributed resources can help enhancing the operation of the power system. Today, the sub-5-minute control of the system is dominated by both primary and secondary frequency control\cite{bergenpower}, \com{where no communication take place among the resources at either level}. Primary frequency control is based on speed droop controllers and the goal is to balance load fluctuations on a short term basis, i.e., within seconds. 
On a slower time scale, but still sub-5-minute, secondary frequency control kicks in to free up the primary control reserves, restore the frequency to 60 Hz, and re-establish the pre-defined tie line flows between control areas. At this level, the communication only takes place between the control center and the resources.

The goal in this paper is to devise a hybrid control structure where a significant part of the control responsibility is transferred to the resources yet with an important contribution of direct coordination among these resources. We particularly focus on secondary frequency control. In current practice, the \emph{independent system operators} (ISOs) first clear the resources via the ancillary service market and then use the  \emph{automatic generation control} (AGC) framework to dispatch the cleared resources according to the pre-defined participation factors. In the wake of the integration of increasing distributed resources \cite{CMKK}, we exploit how to strengthen the communications among the new components to enhance the performance of secondary frequency control. \com{We would like to emphasize that in our paper, we do not require a strongly connected communication network. Instead, the communication among resources can be sparse and limited. And we show that even with such limited communication, the proposed control scheme still stabilizes the system significantly faster than the AGC framework.}

%

To achieve this goal, we utilize the area control error (ACE) signal and propose a secondary frequency control scheme based on a consensus-plus-global-innovation framework \cite{kar2012distributed}. More precisely, in our framework, the participants in real time incrementally update their own control signals in response to the global and via communication exchanged information.
A key objective of the proposed method is to maintain the stability of the system, and at the same time adjust the settings of generators, flexible loads and storages such as to achieve a \emph{cost effective} dispatch, i.e., not just the clearing is optimized but also the dispatch of the cleared resources.

In our earlier work \cite{kar2012distributed2}, a consensus-plus-local-innovation approach is proposed to realize a fully distributed economic dispatch approach, i.e., distributed tertiary control. The procedure is such that the generation units keep exchanging information and updating their solutions until convergence is reached for the predicted load and only the final solution is applied as the optimal dispatch. In the proposed approach here, the swing dynamics of the system make it impossible to wait until convergence. Thus, the control signals are updated and applied immediately after every single iteration. This adds additional complexities to our control design in terms of stability and robustness.

It is worth noting that the application of our approach is not limited to secondary frequency control. Examples of such application abound: e.g., the cost sharing problems in the microgrid, and the electric vehicle charging scheduling, etc.

\vspace{-0.2cm}
\subsection{Related Work}
Frequency stability has been recognized as one of the most important challenges in power systems since the 1970s \cite{converti1976long}. Over the past decade, much effort has been devoted to the design of decentralized robust proportional-integral (PI) controllers to tackle this challenge (see \cite{kumar2005recent} and \cite{pandey2013literature} for two comprehensive surveys). For example, Rerkpreedapong \emph{et al.} combine $H_\infty$ control and genetic algorithm techniques to tune the PI parameters in \cite{rerkpreedapong2003robust}. Bevrani \emph{et al.} introduce a sequential decentralized method to obtain a set of low-order robust controllers in \cite{bevrani2004sequential}. Talaq \emph{et al.} propose an adaptive fuzzy gain scheduling scheme for the conventional PI controller in \cite{talaq1999adaptive}. Chaturvedi \emph{et al.} employ a non-linear neural network controller to perform load frequency control in \cite{chaturvedi1999load}. Li \emph{et al.} propose a modified AGC framework to improve the economic efficiency in \cite{lina2014acc}.
\com{Different from the conventional control perspective with no communications among resources, we utilize the consensus-plus-global-innovation approach and develop a
control framework, which leverages peer-to-peer communication to improve the control performance. We denote such control scheme as distributed communication based control approach.}


Our work also belongs to a growing literature that utilizes the consensus-plus-innovations\footnote{In this paper, since we are using the global information - the ACE signal, we term our approach as consensus-plus-global-innovation.} approach to determine the economic dispatch. For example, Yang \emph{et al.} propose a consensus based distributed framework to perform economic dispatch in \cite{yang2013consensus}. Zhang \emph{et al.} introduce a distributed incremental cost consensus algorithm to solve the economic dispatch problem in \cite{zhang2012convergence}. 
In contrast to this previous work, our approach is significantly more complex in that the outcome after every single iteration is used as a control signal which requires that one has to maintain the system stability while following the dynamically changing demand.

\vspace{-0.2cm}
\subsection{Our Contributions}
In this paper, we seek to answer the following key question: \emph{can we utilize communication among the resources to develop an enhanced secondary frequency control with provable performance?} Towards answering this question, the major contributions of this paper are summarized as follows:

\vspace{0.05cm}
\begin{itemize} \itemsep 2pt
  \item \emph{Consensus-plus-global-innovation Secondary Frequency Control:} Inspired by the economic dispatch problem, we propose a control scheme to perform the secondary frequency control, which satisfies the derived necessary condition to achieve the desirable equilibrium. Each regulation resource updates its local control signal based on the ACE signal received from the control center and information obtained from its communication neighbors.
  \item \emph{System Stability Guarantee:} The proposed control scheme is proven to be a PI controller, which with well developed tuning techniques can stabilize the system. Moreover,  simulation results stress that, thanks to the power of communication, our algorithm can restore the system frequency back to its nominal value significantly faster than the conventional AGC framework does.
  \item \emph{Cost-effectiveness:} Under mild conditions, we prove that the allocations provided by the proposed algorithm are $c\epsilon-$close to the optimally achievable allocations, where $c>0$ is a constant depending only on the cost parameters and the communication topology. This is of particular interest as the proposed scheme does not require the resources to share their cost information with their neighbors nor the control center, which essentially makes it a distributed control scheme.
\end{itemize}

\vspace{0.05cm}

The rest of this paper is organized as follows. We introduce the single-area system model, including the state space model and the time scales of the control scheme in Section \ref{system}. Then, Section~\ref{balancing} demonstrates the consensus-plus-global-innovation control framework. We analyze the properties of our control approach in Section~\ref{analyticalEval}. In Section \ref{gene}, we generalize our control scheme to the multi-area scenario. Simulation results are presented and assessed in Section~\ref{simulation}, confirming the stability and cost effectiveness of
our approach. Section~\ref{conclusions} concludes the paper and points out future works.

\vspace{-0.2cm}
\section{System Model}
\label{system}

\vspace{-0.2cm}
\subsection{Notations}
\begin{itemize}
  \item $\Delta f$: Frequency deviation from nominal value.
  \item $\Delta P_m^i$: Resource $i$'s regulation contribution.
  \item $\Delta P_g^i$: Resource $i$'s governor valve position differential.
  \item $\Delta P_L$: Load deviation from its predicted value.
  \item $H$: Equivalent inertia constant for the single area.
  \item $D$: Equivalent damping coefficient for the single area.
  \item $R_i$: Resource $i$'s droop characteristic.
  \item $u_i$: Secondary frequency control action for generator $i$.
  \item $\Omega=\{1,\cdots,n\}$: The set of all the frequency regulation resources in the system.
  \item \com{$A_G$: Adjacency matrix of the communication graph.
  \item $D_G$: Degree matrix of the communication graph.
  \item $L=D_G-A_G$: The Laplacian matrix of the communication graph.}
\end{itemize}

\vspace{-0.3cm}

\subsection{State Space Model}

\begin{figure}[t]
      \begin{center}
        \includegraphics[width=6.0cm]{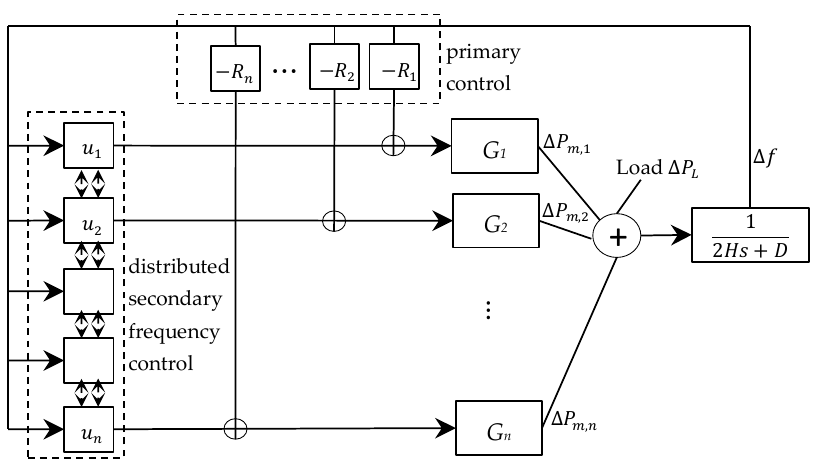}\vspace{-0.3cm}
        \caption{Illustrative single-area control model considered in this paper.}\vspace{-0.8cm}
        \label{fig:system_model}
        \end{center}
\end{figure}

To simplify the analysis and highlight the key properties of our scheme, we first consider the single-area secondary frequency control problem as shown in Fig. \ref{fig:system_model}. The concept differs from the traditional secondary frequency control in that the distributed resources directly communicate with each other and exchange information. We show later that our control scheme can be naturally generalized to the multi-area scenario with tie-line flow constraints.

We assume each frequency regulation resource is a non-reheat steam unit, and use the second order governor-turbine model (in Laplace domain) to characterize it:
\begin{equation}\label{gov-turbine-model}
    M_i(s) = (1+T_g^is)^{-1} (1+T_t^is)^{-1}, \forall i\in\Omega,
\end{equation}
where $T_g^i$ and $T_t^i$ are the governor and turbine time constants for resource $i$, respectively. Thus, in the single area model, we can describe the whole system with the state space model \cite{bergenpower}:
\begin{align}
\Delta \dot{f} \ \ & = -\frac{D}{2H} \Delta f + \frac{1}{2H}\left(\sum_{i\in\Omega} \Delta P_m^i - \Delta P_L\right), \label{freq_swing} \\
\Delta \dot{P}_m^i & = -\frac{1}{T_t^i}(\Delta P_m^i - \Delta P_g^i),\\
\Delta \dot{P}_g^i \ & = -\frac{\Delta f}{T_{g}^i R_i} -\frac{1}{T_g^i}(\Delta P_g^i - u_i).
\end{align}

\vspace{-0.3cm}
\subsection{Time Scales}

Before introducing our discrete time control scheme, we clarify the various system time scales associated with frequency control. The sub-5-minute control time horizon is divided into time slots, $\mathcal{T}=\{1,\cdots,T\}$. Each time slot is of length $\Delta T$. In most major electricity markets in the United States, the current AGC framework sends out the control signal every 4 seconds. Thus, it is natural to set $\Delta T$ to 4 seconds. However, it is worth noting that our control scheme can be applied at even finer time scales. We also make the following assumption so as to formally analyze the performance of our discrete time control scheme using its equivalent continuous time counterpart:

\vspace{0.1cm}
\begin{assumption}
The load deviation $\Delta P_L$ from its prediction only changes at the beginning of each time slot, and remains constant for the rest of the time slot.
\end{assumption}
\vspace{0.1cm}

\begin{figure}[t]
      \begin{center}
        \includegraphics[width=5.2cm]{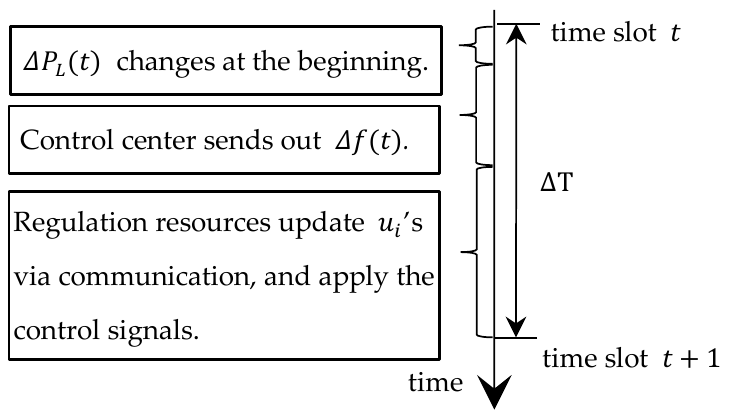}\vspace{-0.3cm}
        \caption{Sequence of actions performed at each time slot.}\vspace{-0.8cm}
        \label{fig:time_scale}
        \end{center}
\end{figure}

The first priority of our control scheme is to update the control signal every $\Delta T$ seconds, and eventually restore the frequency back to the nominal value - 60 Hz, even with $\Delta P_L$ changing at each time slot. We assume, at the beginning of each time slot, the control center sends out the frequency deviation measurement\footnote{\com{In this paper, we employ the standard global frequency deviation measurement. It could be the measurement at certain bus in the system, or the averaged frequency deviation measurement of all buses.}} corresponding to the load change $\Delta P_{L}$ to all the regulation resources. Subsequently, the regulation resources update their own control signals in response to the latest frequency deviation measurement, and also the \emph{information} received from other resources in their respective communication \emph{neighborhoods}. The sequence of actions performed at each time slot is illustrated in Fig. \ref{fig:time_scale}.

\vspace{0.1cm}
\noindent \emph{Remark:} Though Assumption 1 simplifies our later analysis, we would like to point out that the design of our control scheme does not rely on this assumption. It merely utilizes the fact that the current AGC framework sends out the control signal every 4 seconds, which implies that the load variation within the 4-second period is minor. With possibly larger penetration of renewable energy into the system, this assumption may not hold for a period of 4 seconds, but most likely for a period of one second or even smaller. We want to emphasize that both the AGC framework and our proposed control scheme can be naturally generalized to such scenarios if needed.

\vspace{-0.2cm}
\section{Control scheme design}
\label{balancing}

In this section, we study the (controlled) swing dynamics and obtain necessary conditions that any control scheme needs to satisfy in order to achieve desired equilibrium behavior, i.e., equilibrium points at which the frequency deviation $\Delta f$ from the nominal is zero.  Subsequently, we investigate the secondary frequency control objective from an economic dispatch perspective. Motivated by these two control objectives, i.e., desired equilibrium behavior and cost-effectiveness of dispatch, we introduce our control scheme, which leverages peer-to-peer communication.

\vspace{-0.2cm}
\subsection{Necessary Conditions  to Achieve Desired Equilibrium}

Any equilibrium point of the  controlled swing dynamic equations should satisfy that \com{ $\Delta \dot{f}, \ \Delta \dot{P}_m^i,$ and $\Delta \dot{P}_g^i$ are all zero.}
We want to design a control scheme, i.e., determining the $u_i$'s, such that, we have $\Delta f=0$ at equilibrium. Thus, necessary conditions to achieve such an equilibrium are
\begin{equation}\label{ness1}
    \begin{cases}
    \sum_{i\in\Omega} \Delta P_m^i = \Delta P_L,\\
    \Delta P_m^i = \Delta P_g^i = u_i.
    \end{cases}
\end{equation}
In particular, any control scheme that wants to achieve an equilibrium with $\Delta f =0$ should at least ensure
\begin{equation}\label{ness}
    \textstyle \sum_{i\in\Omega} u_i = \Delta P_L.
\end{equation}
This is a well known fact in power system operation and should be enforced if a given frequency control scheme is to restore the system frequency back to its nominal value.

\vspace{-0.2cm}
\subsection{Economic Dispatch Inspiration}

It is readily seen that there could be many control schemes satisfying the necessary condition (\ref{ness}), but we want to design a control scheme that also achieves cost effectiveness in the sense of the cost of energy dispatch.

Mathematically, if we consider a quadratic cost function for each regulation resource $i$, from an economic perspective, ideally, we would want, by the end of time slot $t$, $P_m^i(t+1)$'s that minimize the dispatch cost:
\begin{align}
   & \text{\textbf{minimize}} \ \ \textstyle\sum_{i\in\Omega} \left ( a_i (\Delta P_m^i)^2 + {b}_i \Delta P_m^i + c_i \right ),\label{p_quadcost}\\
   & \text{\textbf{subject to }} \textstyle\sum_{i\in\Omega} \Delta P_m^i = \Delta P_L(t+1). \label{p_balance}
\end{align}

Note, in the above we ignore the ramping constraints of the network resources which decouples the temporal constraints, thereby reducing the overall dispatch cost minimization over the control horizon $\mathcal{T}$ to solving a collection of temporally decoupled economic dispatch problems (\ref{p_quadcost})-(\ref{p_balance}) at each time slot $t$. In fact, as long as $|\Delta P_L(t)-\Delta P_L(t-1)|$ is appropriately bounded throughout $\mathcal{T}$, this relaxation is tight. We provide this justification in Appendix \ref{app1}.

Further, note that since the regulation resources are selected by co-optimizating with the energy bids, taking the capacity payment (used to compensate the lost opportunity cost in the energy market\cite{federal2011frequency}) into account, the initial marginal costs (when $\Delta P_m^i=0$), i.e., $b_i$'s, of all the participants to provide secondary frequency control should be almost identical. Therefore, we make the following assumption to simplify the analysis:

\vspace{0.1cm}
\begin{assumption}
All the $b_i$'s are identical, i.e., $b_i=b, \forall i\in\Omega$.
\end{assumption}
\vspace{0.1cm}

We know at the equilibrium point, we have $u_i=\Delta P_m^i$, for all $i\in\Omega$. Combining this with (\ref{p_quadcost})-(\ref{p_balance}), ideally, from a cost of dispatch viewpoint, we would like to design a control scheme, such that $u_i(t+1)$'s are the solution to
\begin{align}
   & \text{\textbf{minimize}} \ \ \textstyle\sum_{i\in\Omega} \left ( a_i  u_i^2 + {b}  u_i + c_i \right ),\label{p_quadcost_u_c}\\
   & \text{\textbf{subject to }} \textstyle\sum_{i\in\Omega}  u_i = \Delta P_L(t+1). \label{p_balance_u_c}
\end{align}
Note that, with identical $b_i$'s, constraint (\ref{p_balance_u_c}) guarantees $\sum_{i\in\Omega}{b}  u_i=b\Delta P_L(t+1)$, which is a constant in the optimization problem. Hence, by reducing all the constant terms, optimization problem (\ref{p_quadcost_u_c})-(\ref{p_balance_u_c}) reduces to
\begin{align}
   & \text{\textbf{minimize}} \ \ \textstyle\sum_{i\in\Omega}  a_i  u_i^2 ,\label{p_quadcost_u}\\
   & \text{\textbf{subject to }} \textstyle\sum_{i\in\Omega}  u_i = \Delta P_L(t+1). \label{p_balance_u}
\end{align}

Denoting all the control signals $\{u_i(t+1),\forall i\in\mathcal{N}\}$ by $\boldsymbol{u}(t+1)$, and the Lagrangian multiplier associated with (\ref{p_balance_u}) by $\lambda^{t+1}$, the Lagrangian function $\mathcal{L}(\boldsymbol{u}(t+1),\lambda^{t+1})$ for problem (\ref{p_quadcost_u})-(\ref{p_balance_u}) at a given time $t$ is given by
\begin{equation}\label{lag}
\begin{aligned}
   & &  &\mathcal{L}(\boldsymbol{u}(t+1),\lambda^{t+1}) \\
   & & = &\! \sum_{i\in\Omega}\! a_i  u_i^2(t\!+\!1)\!-\! \lambda^{t+1}\!\left ( \sum_{i\in\Omega}\!  u_i(t\!+\!1) \!-\! \Delta P_L(t\!+\!1)\!\! \right)\!.
    \end{aligned}
\end{equation}
The first order optimality conditions are therefore given by
\begin{align}
    & \frac{\partial \mathcal{L}}{\partial u_i(t+1)} = 2a_i u_i(t+1) - \lambda^{t+1}=0, \ \forall i\in\Omega,\label{sec3:econ_optimal} \\
    & \frac{\partial \mathcal{L}}{\partial \lambda^{t+1}} = \sum_{i\in\Omega} u_i(t+1) - \Delta P_L(t+1)=0. \label{sec3:balance}
\end{align}
If we denote each regulation resource $i$'s marginal cost at time $t$ by $\lambda_i^t$, (\ref{sec3:econ_optimal}) requires
\begin{equation}\label{econ_con}
    \lambda_i^{t+1} = 2a_i u_i(t+1)= \lambda^{t+1}, \ \forall i\in\Omega.
\end{equation}
Hence, (\ref{sec3:econ_optimal})-(\ref{econ_con}) reflect that the marginal costs for all entities have to be equal in the optimal solution and the resulting provision of power needs to fulfill the power balance.

At a given time $t\in\mathcal{T}$ if each resource $i$ has access to the Lagrange multiplier variable $\lambda^t$, which may be interpreted as a differential generation price, it may set its control signal according to (\ref{sec3:econ_optimal}) so that the system achieves the power balance in the most economic way. However, as (\ref{sec3:balance}) suggests, the quantity $\lambda^t$ depends on private information such as the cost characteristics of all the entities and global information - the instantaneous system net load deviation $\Delta P_L(t)$. Thus, we naturally ask the following question: \emph{is it possible to design a distributed control scheme which does not require all the private information yet can still achieve provable performance?} This motivates us to propose the consensus-plus-global-innovation control scheme.

\vspace{-0.2cm}
\subsection{Consensus + Global Innovation Control Design}



We propose a distributed real-time approach in which participating entities, through neighborhood communication and global information processing, continuously update their control signals to track the optimal power allocation closely.
Note that we denote our control scheme as \textit{distributed} for the following reasons: cost parameters of resource $i$ are only known to resource $i$ and communication is used to find an agreement with neighboring resources.
In our control scheme,  each resource $i\in\Omega$ maintains and updates a local copy of the variable $\lambda_i^t$. The updates are defined as
\begin{align}
 & &  &  \tilde{\lambda}_i^{t+1} = {\lambda}_i^t -\overset{\text{neighborhood consensus}}{ \overbrace{2{a}_i \beta \textstyle\sum_{l\in\Omega_i} ({\lambda}_i^t -{\lambda}_l^t)}}\label{sec3:lambda_update}\\
& & & \ \ \ \ \ \ \ \ \ \ +\!\underset{\text{global innovation}}{ \underbrace{2{a}_i n^{-1}(\Delta P_L(t\!+\!1)\! -\!  \textstyle\sum_{j\in\Omega}\Delta{P}_m^j(t))}},\nonumber \\
 & & & u_i(t+1) = (2a_i)^{-1} \tilde{\lambda}_i^{t+1}, \\
 & &  & \lambda_i^{t+1} =  2a_i \Delta {P}_m^i(t+1), \label{sec3:lambda_reupdate}
\end{align}
where $\beta>0$ is a tuning parameter; \com{$\tilde{\lambda}_i^t$ is the estimate of marginal cost $\lambda_i^t$ for resource $i$ at time slot $t$;} $\Omega_i$ denotes the set of participant $i$'s neighbors in the communication network\footnote{We assume that an inter-resource communication network is pre-defined. Moreover, note that this communication network could be different from the physical power system topology and possibly much sparser.}.

Intuitively, note that the neighborhood consensus term in the update rule (\ref{sec3:lambda_update}) seeks to enforce an agreement between the marginal price variables $\lambda_{i}^{t}$'s so as to optimize the dispatch cost (see (\ref{econ_con})); whereas, the innovation term seeks to enforce demand-supply balance which is also necessary to drive $\Delta f$ to zero.
This update would require that entity $i$ has access to  $\Delta P_L(t)$ and all $\Delta P_m^i(t)$'s, which constitutes global as opposed to local information. In order to realize the (global) innovation term in (\ref{sec3:lambda_update}) using local information, we use the fact that
\begin{equation}\label{Plest_c}
   \Delta P_L = \textstyle \sum_{i\in\Omega} \Delta P_m^i -D \Delta f- 2H \Delta \dot{f},
\end{equation}
which follows from the swing dynamics (\ref{freq_swing}). The discrete time approximation of $\Delta P_L$ therefore is
\begin{equation}\label{PLest}
\begin{aligned}
   \Delta P_L (t\!+\!1) =& \textstyle \sum_{i\in\Omega} \Delta P_m^i(t) -D \Delta f(t)\\
    & \ \ - 2H(\Delta T)^{-1} (\Delta f(t\!+\!1)\!-\!\Delta f(t)).
   \end{aligned}
\end{equation}
Substituting (\ref{PLest}) into the global innovation term in (\ref{sec3:lambda_update}) yields
\begin{equation}\label{sec3:lambda_update_renew}
    \begin{aligned}
   & \tilde{\lambda}_i^{t+1} = {\lambda}_i^t -\overset{\text{consensus potential}}{ \overbrace{2{a}_i \beta \textstyle\sum_{l\in\Omega_i} ({\lambda}_i^t -{\lambda}_l^t)}}\\
   & \ \ \ \ \ \ \ \ \ \underset{\text{innovation potential}}{ \underbrace{-\frac{2{a}_i}{n}D \Delta f(t)-\frac{4{a}_iH}{ n\Delta T}(\Delta f(t\!+\!1)\!-\!\Delta f(t))}}.
    \end{aligned}
\end{equation}

\section{Analytical Performance Evaluation}
\label{analyticalEval}

Given this distributed control scheme, we analyze its key properties - stability, equilibrium behavior, and dispatch cost effectiveness in this section.

\vspace{-0.2cm}
\subsection{Stability}

For each $u_i(t+1)$, we can rewrite the control law as
\begin{equation}\label{control_law_ui}
\begin{aligned}
    u_i(t+1) = & \Delta P_m^i(t)-\beta L^i\Lambda^{-1} \Delta \boldsymbol{P}_m(t) \\
    & + n^{-1}\left ({\Delta P_L(t\!+\!1)}\!-\! \boldsymbol{1}^T \boldsymbol{P}_m(t)\right ),
    \end{aligned}
\end{equation}
where $\boldsymbol{P}_m(t)=[\Delta P_m^1(t),\cdots,\Delta P_m^n(t)]^T$; $L^i$ is the $i^{\text{th}}$ row of the communication network's Laplacian matrix $L$; 
and
\begin{equation}
   \Lambda = \text{diag} \{(2a_1)^{-1}, \cdots, (2a_n)^{-1}\},
\end{equation}
where diag$\{\cdot \}$ denotes the diagonal matrix.

Using the Euler forward emulation and (\ref{PLest}), we have the Laplace transform for $u_i$ as
\begin{equation}\label{control_ui_lap}
\begin{aligned}
    u_i(s)(1+\Delta Ts) = &  \Delta P_m^i(s)-\beta L^i\Lambda \Delta \boldsymbol{P}_m(s)\\
 &    - n^{-1}(2H\Delta f(s)s+D\Delta f(s)).
 \end{aligned}
\end{equation}

Note that, according to the second order model (\ref{gov-turbine-model}), we can further approximate $\Delta P_m^i(s)$ by
\begin{equation}\label{Pmfurtherappro}
\begin{aligned}
    \Delta P_m^i(s) &  = u_i(s) (1+T_g^is)^{-1}(1+T_t^is)^{-1} \\
    & \approx u_i(s)(1-(T_g^i+T_t^i)s).
    \end{aligned}
\end{equation}
Combining (\ref{Pmfurtherappro}) with (\ref{control_ui_lap}) yields
\begin{equation}\label{control_ui_lap_final}
    u_i(s) = - \frac{2H}{nT_u^i}\Delta f(s) - \frac{\beta nL^i\Lambda \Delta \boldsymbol{P}_m(s) + D\Delta f(s) }{nT_u^is},
\end{equation}
where $T_u^i = \Delta T +T_g^i+T_t^i$.

That is, our control may indeed be approximated as a PI controller.
Given the representation (\ref{control_ui_lap_final}), we can use standard results on control, see for example \cite{cao1998static}, to design the various algorithm parameters in order to appropriately shape the closed loop frequency response of the system and, in particular, to ensure stability. While we do not pursue this in detail in this paper, some generic guidelines for tuning the algorithm parameters to achieve stable closed loop behavior are provided in Section \ref{simulation}.

\vspace{0.1cm}
\noindent{\emph{Remark}:} Compared with the conventional AGC framework, our distributed control scheme makes use of more state information (both $\Delta f$ amd $\Delta P_m^i$'s) via communication. This grants us additional flexibility in terms of designing the (PI) controller. As demonstrated by the simulation results (see Section \ref{simulation}), the power of communication enables our control scheme to stabilize the system significantly faster than the AGC framework.

\vspace{-0.2cm}
\subsection{Equilibrium Behavior}

We have mentioned that our control scheme satisfies the necessary conditions to achieve the desired equilibrium behavior. Though it is the standard requirement for power system control, we formally illustrate this property in Theorem \ref{thm:realtimebalance}.

\vspace{0.1cm}
\begin{theorem}\label{thm:realtimebalance}
The proposed distributed control scheme satisfies the necessary condition (\ref{ness}).
\end{theorem}
\vspace{0.1cm}

\noindent \emph{Proof:} At each time $t\in\mathcal{T}$, a step load change may happen at the beginning of the time slot, and according to (\ref{sec3:lambda_update})-(\ref{sec3:lambda_reupdate}), we have
\begin{equation}\label{ind_1}
   \textstyle \sum_{i\in\Omega} u_i(t+1) = \sum_{i\in\Omega}{(2a_i)^{-1}} {\tilde{\lambda}_i^{t+1}},
\end{equation}
and
\begin{equation}\label{sec3:induction_P}
    \begin{aligned}
    & &  \sum_{i\in\Omega}\!{(2a_i)^{-1}} \tilde{\lambda}_i^{t+1} \!=\!&\sum_{i\in\Omega}\!\! {(2a_i)^{-1}} \tilde{\lambda}_i^{t}\! -\!\beta\! \sum_{i\in\Omega} \! \sum_{l\in\Omega_i}\! (\lambda_i^t\!-\!\lambda_l^t) \\
    & & & +\! {\Delta P_L(t+1)} - \sum_{i\in\Omega} \! \Delta{P}_m^i(t).
    \end{aligned}
\end{equation}
Note the following facts hold:
\begin{equation}\label{bala_eq}
   \textstyle \sum_{i\in\Omega} \sum_{l\in\Omega_i}\! (\lambda_i^t\!-\!\lambda_l^t) = 0, \forall t\in\mathcal{T};
\end{equation}
\begin{equation}\label{bala_eq2}
\textstyle \sum_{i\in\Omega} (2a_i)^{-1}{\lambda_i^t} = \sum_{i\in\Omega} \! \Delta{P}_m^i(t).
\end{equation}
Combining (\ref{bala_eq}) and (\ref{bala_eq2}) with (\ref{sec3:induction_P}), we establish that
\begin{equation}\label{conclu_thm2}
 \textstyle \sum_{i\in\Omega}\! u_i(t\!+\!1) = \Delta P_L(t+1).
\end{equation}
$\hfill \blacksquare$
\vspace{-0.4cm}

\subsection{Cost Effectiveness}

We analyze the cost effectiveness of the proposed distributed control scheme under the following assumption.

\vspace{0.1cm}
\begin{assumption}\label{zero-time}
If the time constants $T_g^i$'s and $T_t^i$'s in the second order model (\ref{gov-turbine-model}) are all zero, we have for all $t\in\mathcal{T}$
\begin{equation}\label{Pmappro_1}
    \Delta P_m^i(t) = u_i(t), \forall i\in\Omega.
\end{equation}
\end{assumption}

With increasing storage systems, and demand response resources participating in the secondary frequency control market, more resources enjoy very small time constants. This control scheme can thus guarantee near real time demand-supply balance. We rely on this assumption to prove the cost effectiveness of our control scheme, but for the simulations we still use reasonable time constants for non-reheat steam turbines. This assumption basically corresponds to saying that the control settings for the generators are optimal which after the physical delay caused by the governor and the turbine will be seen at the output after this delay.

To prove cost effectiveness, we also need the communication graph to be well connected. In particular, the second largest eigenvalue $1-\rho$ of matrix $I-\beta\Lambda^{-1} L$ need satisfy:
\begin{equation}\label{thm2-assum}
    (1-\rho)\textstyle\max_{i\in\mathcal{N}} (a_i)^{\frac{1}{2}}\textstyle \max_{i\in\mathcal{N}} (a_i)^{-\frac{1}{2}} < 1.
\end{equation}

\vspace{0.1cm}\com{
\noindent \emph{Remark:} Condition (\ref{thm2-assum}) guarantees sufficient speed to implicitly spread the cost parameters via consensus. The Lapalacian matrix has a trivial eigenvalue of zero, corresponding to the largest eigenvalue, 1, of the matrix $I-\beta\Lambda^{-1} L$. The second largest eigenvalue of this matrix measures the connectivity of the communication network. The more strongly connected the network, the smaller the second largest eigenvalue ($1-\rho$). Note that, if $a_i$'s are all the same, there is no need to spread the cost parameters through consensus. We merely require $1-\rho$ is strictly less than 1, which can be easily achieved as long as the network is connected. On the other hand, if the $a_i$'s vary a lot, we need a reasonably connected network to well spread such information indirectly via the $\lambda_i^t$'s. }

\vspace{0.1cm}
\begin{theorem}\label{thm:costeffective}
If the communication graph is well connected such that (\ref{thm2-assum}) holds, $\|\Delta P_L(t+1)- \Delta P_L(t)\| < \epsilon$, $\forall t\in \mathcal{T}$, and Assumption \ref{zero-time} holds, then there exists a constant $c > 0$ (depending only on the cost parameters and
the communication topology) such that  control signals ($u_i$'s) are $c\epsilon$-close to their optimal values given by (\ref{sec3:econ_optimal})-(\ref{sec3:balance}). More precisely,
\begin{align}\label{sec3:cost_goal}
&   \|\lambda_i^t - \lambda^t\| \le c\epsilon, \ \forall t\in\mathcal{T}, i\in\Omega,\\
  &   \|u_i(t) - u^*_i(t)\| \le c\epsilon, \ \forall t\in\mathcal{T}, i\in\Omega,
  \end{align}
where $u^*_i(t)$ denotes the solution to (\ref{sec3:econ_optimal})-(\ref{sec3:balance}).
\end{theorem}
\vspace{0.1cm}

\noindent \emph{Proof:} Using Assumption \ref{zero-time}, we know that $\lambda^t = \tilde{\lambda}^t$, for all $t\in\mathcal{T}$. Thus, the marginal cost updates in (\ref{sec3:lambda_update}) may be written in vector notation as
\begin{equation}\label{lambda_vec}
   \Lambda {{\boldsymbol{\lambda}}}_{t+1} = (\Lambda -\beta L ) {\boldsymbol{\lambda}}_{t} + n^{-1}(\Delta P_L(t+1) -\Delta P_L(t)) \boldsymbol{1}.
\end{equation}
Since $\boldsymbol{1}^T L = 0^T$, we have
\begin{equation}\label{eq0}
    \boldsymbol{1}^T \Lambda {\boldsymbol{\lambda}}_{t+1} =  \boldsymbol{1}^T \Lambda {\boldsymbol{\lambda}}_{t} + (\Delta{P}_L(t+1)-\Delta{P}_L(t)).
\end{equation}
Assumption \ref{zero-time} also leads to
\begin{equation}\label{eq2}
   \boldsymbol{1}^T \Lambda {\boldsymbol{\lambda}}_{t} = \boldsymbol{1}^T \Lambda  \boldsymbol{1} {\lambda}^t.
\end{equation}
Combining (\ref{eq2}) with (\ref{eq0}) yields
\begin{equation}\label{eq3}
    \boldsymbol{1}^T \Lambda  \boldsymbol{1} {\lambda}^{t+1} = \boldsymbol{1}^T \Lambda  \boldsymbol{1} {\lambda}^t + (\Delta{P}_L(t+1)-\Delta P_L(t)).
\end{equation}
By dividing $\boldsymbol{1}^T \Lambda  \boldsymbol{1} $ at both sides, we have
\begin{equation}\label{eq4}
     {\lambda}^{t+1} = {\lambda}^{t} + (\boldsymbol{1}^T \Lambda  \boldsymbol{1})^{-1}(\Delta{P}_L(t+1)-\Delta P_L(t)).
\end{equation}
\com{Therefore,
\begin{equation}\label{main_eq}
    \begin{aligned}
  \!\!\!\!  & &  & {\boldsymbol{\lambda}}_{t+1} -  {\lambda}^{t+1} \boldsymbol{1} \\
  \!\!\!\!  & & = & (I\!-\!\beta \Lambda^{-1} L ) {\boldsymbol{\lambda}}_{t}\! +\! (n\Lambda)^{-1}(\Delta P_L\!(t\!+\!1)\!-\! \Delta{P}_{L}\!(t)) \boldsymbol{1}\! \\
  \!\!\!\!   & & &  - \! {\lambda}^t  \boldsymbol{1} \!-  \! ( \boldsymbol{1}^T \Lambda  \boldsymbol{1})^{-1}(\Delta{P}_L(t+1)-\Delta P_L(t))  \boldsymbol{1} \\
  \!\!\!\!  & & = & (I\!-\! \beta \Lambda^{-1} L ) ({\boldsymbol{\lambda}}_{t} \!- \!{\lambda}^t\boldsymbol{1}  ) \\
  \!\!\!\!   & & & + ((n\Lambda)^{-1} \!-  \! ( \boldsymbol{1}^T \Lambda  \boldsymbol{1})^{-1})(\Delta{P}_L(t+1)-\Delta P_L(t))  \boldsymbol{1} \\
  \!\!\!\!  & & \le & (I\!-\! \beta \Lambda^{-1} L ) ({\boldsymbol{\lambda}}_{t} \!- \!{\lambda}^t\boldsymbol{1}  ) \! +\! \|(n\Lambda)^{-1}-( \boldsymbol{1}\!^T\! \Lambda  \boldsymbol{1})\!^{-1}\|\epsilon\boldsymbol{1}.\\
    \end{aligned}
\end{equation}
Note that in the second equation of (\ref{main_eq}), we add a zero term $\beta \Lambda^{-1} L {\lambda}^t\boldsymbol{1}$.}
By contradiction, we can show that the largest eigenvalue of the matrix $I-\beta \Lambda^{-1} L $ is $1$.
Note that,
\begin{equation}\label{fix_1}
I- \beta \Lambda^{-1}  L = \Lambda^{-\frac{1}{2}} (I-   \beta \Lambda^{-\frac{1}{2}} L \Lambda^{-\frac{1}{2}})  \Lambda^{\frac{1}{2}},
\end{equation}
$I\!-\! \beta \Lambda^{-1}   L$ and $I\!-\!   \beta \Lambda^{-\frac{1}{2}} L \Lambda^{-\frac{1}{2}}$ are similar matrices with the same eigenvalues. Since the Laplacian $L$ corresponds to a connected network, one is a simple eigenvalue of the matrix $I-\beta\Lambda^{-1/2}L\Lambda^{-1/2}$ with corresponding eigenvector $\Lambda^{\frac{1}{2}} \boldsymbol{1}$. By (\ref{eq2}), we have
\begin{equation}\label{orth}
    (\Lambda^{\frac{1}{2}} \boldsymbol{1})^T \Lambda^{\frac{1}{2}}({\boldsymbol{\lambda}}_{t} -  {\lambda}^t\boldsymbol{1} ) = 1^T \Lambda ({\boldsymbol{\lambda}}_{t} -  {\lambda}^t\boldsymbol{1} ) = 0.
\end{equation}
Thus, $ \Lambda^{\frac{1}{2}}({\boldsymbol{\lambda}}_{t} -  {\lambda}^t\boldsymbol{1})$ is orthogonal to the eigenvector corresponding to the eigenvalue 1. Furthermore, the matrix $\Lambda^{-1/2}L\Lambda^{-1/2}$ is positive semidefinite and hence, by taking $\beta$ to be small enough, all other eigenvalues of $I-\beta\Lambda^{-1/2}L\Lambda^{-1/2}$ can be guaranteed to lie in the interval [0,1).  Since $I\!-\!   \beta \Lambda^{-\frac{1}{2}} L \Lambda^{-\frac{1}{2}}$ is a symmetric matrix, the induced 2-norm coincides with the spectral radius, thus
\begin{equation}\label{fix_2}
\begin{aligned}
& & & \|(I-  \Lambda^{-1}  \beta L ) ({\boldsymbol{\lambda}}_{t} - {\lambda}^t\boldsymbol{1})  \|  \\
& & = & \|\Lambda^{-\frac{1}{2}}(I-   \beta \Lambda^{-\frac{1}{2}} L \Lambda^{-\frac{1}{2}}) \Lambda^{\frac{1}{2}} ({\boldsymbol{\lambda}}_{t} - {\lambda}^*(t)\boldsymbol{1} ) \|\\
& & \le & \|\Lambda^{-\frac{1}{2}}\| \|I-   \beta \Lambda^{-\frac{1}{2}} L \Lambda^{-\frac{1}{2}}\| \| \Lambda^{\frac{1}{2}}\| \|{\boldsymbol{\lambda}}_{t}- {\lambda}^t\boldsymbol{1}  \|\\
& & \le &  (1-\rho) \textstyle\max_{i\in\mathcal{N}} a_i^{\frac{1}{2}} \textstyle\max_{i\in\mathcal{N}} (a_i)^{-\frac{1}{2}} \|{\boldsymbol{\lambda}}_{t} - {\lambda}^t \boldsymbol{1} \|,
\end{aligned}
\end{equation}
where $1-\rho\in [0,1)$ is  the second largest eigenvalue of $I-\beta \Lambda^{-1} L $. As long as (\ref{thm2-assum}) holds, standard algebraic
manipulations now lead to the desired conclusion.
$\hfill \blacksquare$

Theorem \ref{thm:costeffective} guarantees that, as long as the communication topology is well connected \com{(satisfying (\ref{thm2-assum}))}, the deviations of the control signals $u_i$'s from their optimal settings $u_i^*(t)$'s (in the sense of the economic dispatch objective) at any time $t\in\mathcal{T}$ are upper-bounded by a quantity proportional to the maximum rate of change of demand (i.e., $\epsilon$).
Since the rate of change of demand is usually fairly small, this theorem implies that the allocations stay close to their optimal at all times and hence are cost-effective.

\vspace{-0.2cm}
\section{Multi-area Generalization}
\label{gene}

Though our analysis is based on the single area model, where the ACE signal is exactly the frequency deviation, it is not hard to see that by replacing the current control signal $\Delta f(t)$ with the ACE signal, we can directly apply the proposed distributed secondary frequency control scheme to the multi-area model. More specifically, the frequency swing dynamics in area $j$ becomes
\begin{equation}\label{freq_sw_multi}
    \Delta \dot{f}_j \!=\! -\frac{1}{2H_j}\!\! \left(\!\!D_j\Delta f_j \!+\!\Delta P_{tie,j}\! -\!\!\!\sum_{i\in\Omega_j} \Delta P_{m,j}^i \!+\!\Delta P_{L,j}\!\!\right)\!\!,\!
\end{equation}
where $\Delta P_{tie,j}$ denotes the net tie-line flow out of area $j$. Thus, the global innovation term in the updating rule (\ref{sec3:lambda_update}) can be estimated by
\begin{equation}\label{P_lest_multi}
\begin{aligned}
    & \Delta P_{L,j}(t+1) -\textstyle\sum_{i\in\Omega_j} \Delta P_{m,j}^i(t)\\
     = &-\!2H_j\frac{\Delta {f}_j(t\!+\!1)\!-\!\Delta f_j(t)}{\Delta T}\! -\!D_j\Delta f_j(t)\! -\!\Delta P_{tie,j}(t).\!
   \end{aligned}
\end{equation}
In other words, each area will update its own innovation term in (\ref{sec3:lambda_update}) according to (\ref{P_lest_multi}). Based on this information and the communication with the neighbors in the same control area, the regulation resources update the control signals.
The performance analysis is essentially the same as the single-area scenario, which we omit due to space constraints. We confirm the stability of our distributed control framework in the multi-area scenario with simulations in the next section.

\vspace{-0.2cm}
\section{Simulation Results}
\label{simulation}

In this section, we first carry out the simulation for the proposed control scheme for a single-area system, and then test our control scheme in the multi-area scenario and compare the performance to the performance of the conventional AGC scheme.

\vspace{-0.2cm}
\subsection{Single-area Scenario}

In the single-area five-participant system, the equivalent inertia constant $H$ is assumed to be $0.0833$ pu s; and the equivalent damping coefficient $D$ is set to $0.0084$ pu/Hz. In the updating rules, we choose $\beta=0.003$. For the AGC simulations, we use uniform participation factors, and select the coefficients for its PI controller such that it leads to the best achievable performance. The droop characteristics $R_i$'s, the governor  time constants $T_g^i$'s, and the turbine time constants $T_t^i$'s are generated uniformly at random from [2,3] Hz/pu, [0.05,0.06] s, and [0.3,0.5] s, respectively. The communication network is a 2-nearest neighbor network.

\begin{figure}[t]
      \begin{center}
        \includegraphics[width=5.4cm]{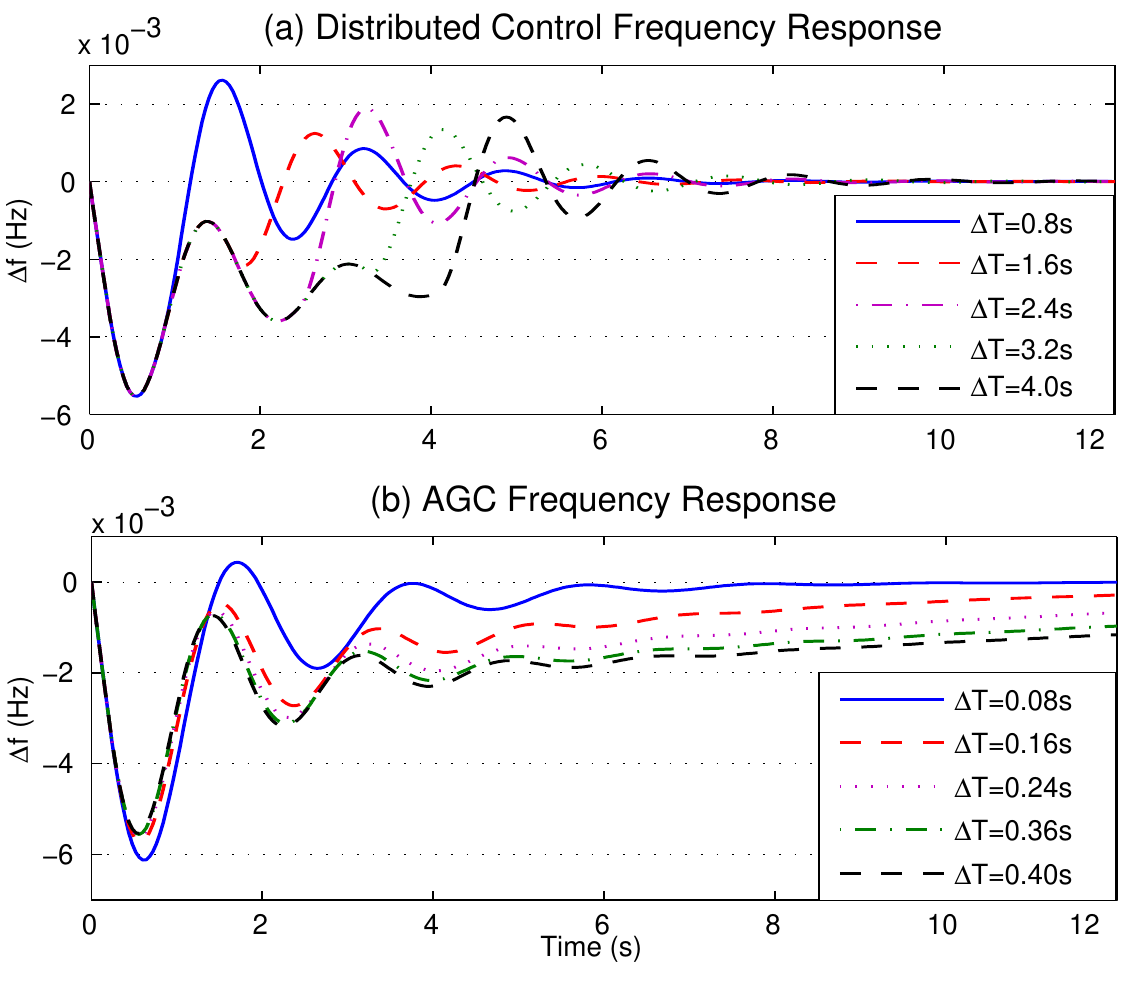}\vspace{-0.4cm}
        \caption{$\Delta T$ impact of the stability of two systems.}\vspace{-0.8cm}
        \label{fig:DeltaTImpact}
        \end{center}
\end{figure}

We first study the impact of the control interval length $\Delta T$ on the stability of both our distributed control scheme and the AGC framework. As shown in Fig. \ref{fig:DeltaTImpact}, when facing a step load increase of 0.005 pu, our distributed control scheme outperforms the AGC framework with varying $\Delta T$: The AGC framework takes significantly longer than our scheme to restore the frequency back to 60 Hz with larger $\Delta T$.

More specifically, when $\Delta T=4$s, our control scheme can stabilize the system within 8s - that is, with two control signal updates. On the other hand, even when $\Delta T=0.16$s, the current AGC system needs more than 12 seconds to stabilize the system. We ascribe our control system's advantage on fast stabilizing the system to its utilization of communications and more state information (both $\Delta f$ and $\Delta P_m^i$'s) while the AGC framework only relies on the frequency information.

\begin{figure}[t]
      \begin{center}
        \includegraphics[width=5.45cm]{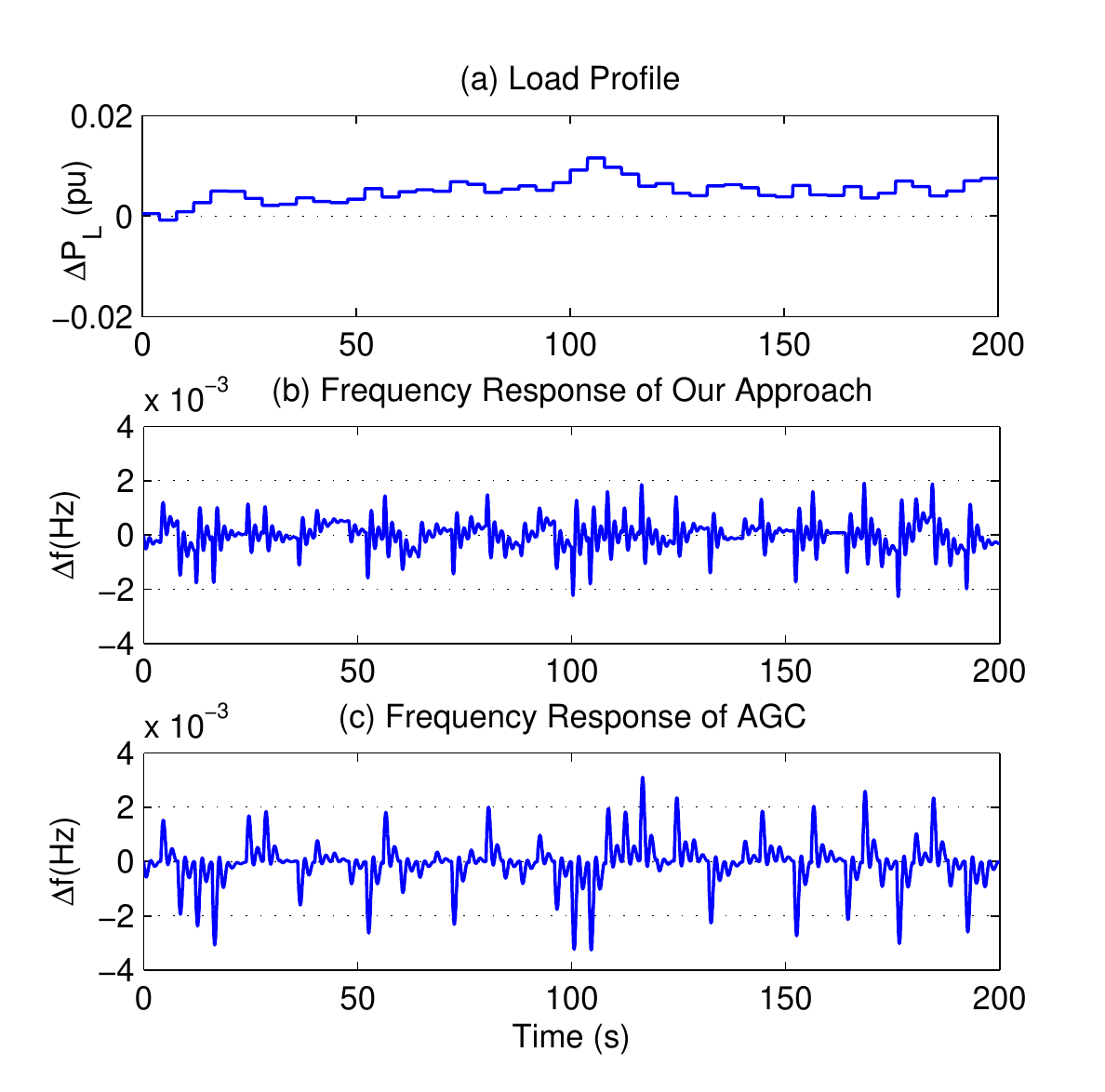}\vspace{-0.3cm}
        \caption{Frequency responses with continuously changing load: $\Delta T$=0.08s.}\vspace{-0.4cm}
        \label{fig:const_PL_t008}
        \end{center}
\end{figure}

\begin{figure}[t]
      \begin{center}
        \includegraphics[width=5.4cm]{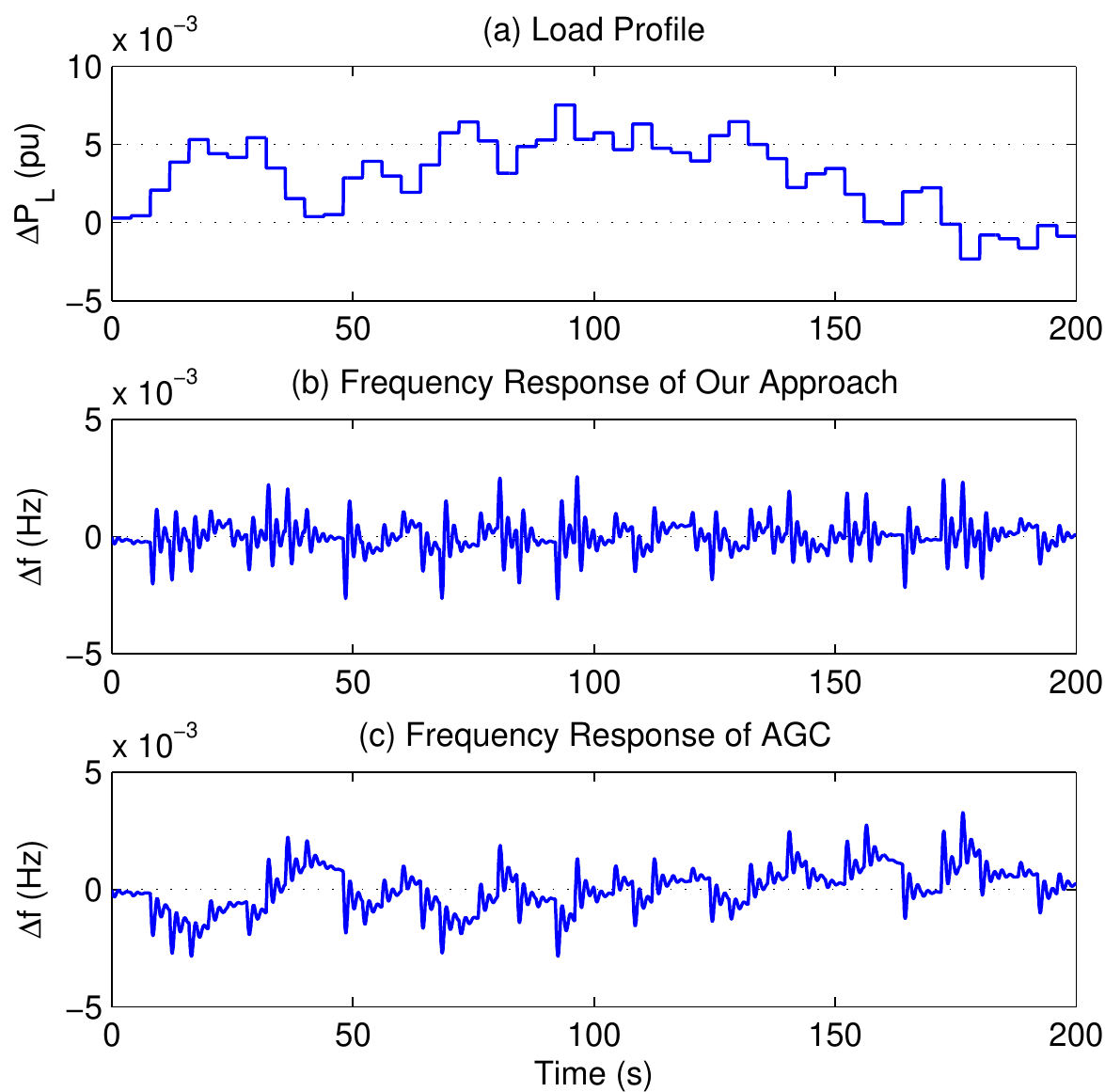}\vspace{-0.3cm}
        \caption{Frequency responses with continuously changing load: $\Delta T$=0.4s.}\vspace{-0.8cm}
        \label{fig:const_PL_t04}
        \end{center}
\end{figure}

Next, by choosing common $\Delta T$ as 0.08s and 0.4s, we compare the frequency responses of the two systems when the load constantly changes, as demonstrated in Fig. \ref{fig:const_PL_t008} and \ref{fig:const_PL_t04}, respectively. Here, we also use Assumption 1, that the load only changes every four seconds. When $\Delta T$ is set to be 0.08s, the frequency responses of the two control schemes both work reasonably well. However, when $\Delta T$ is 0.4s, the AGC framework needs significantly longer time to restore the frequency. Therefore, as shown in Fig. \ref{fig:const_PL_t04}(c), the frequency response of the AGC framework in this case constantly deviates from the nominal value. It is interesting to note that, our control scheme in this case works as well as the case when $\Delta T$ is 0.08s.

\begin{figure}[t]
      \begin{center}
        \includegraphics[width=5.5cm]{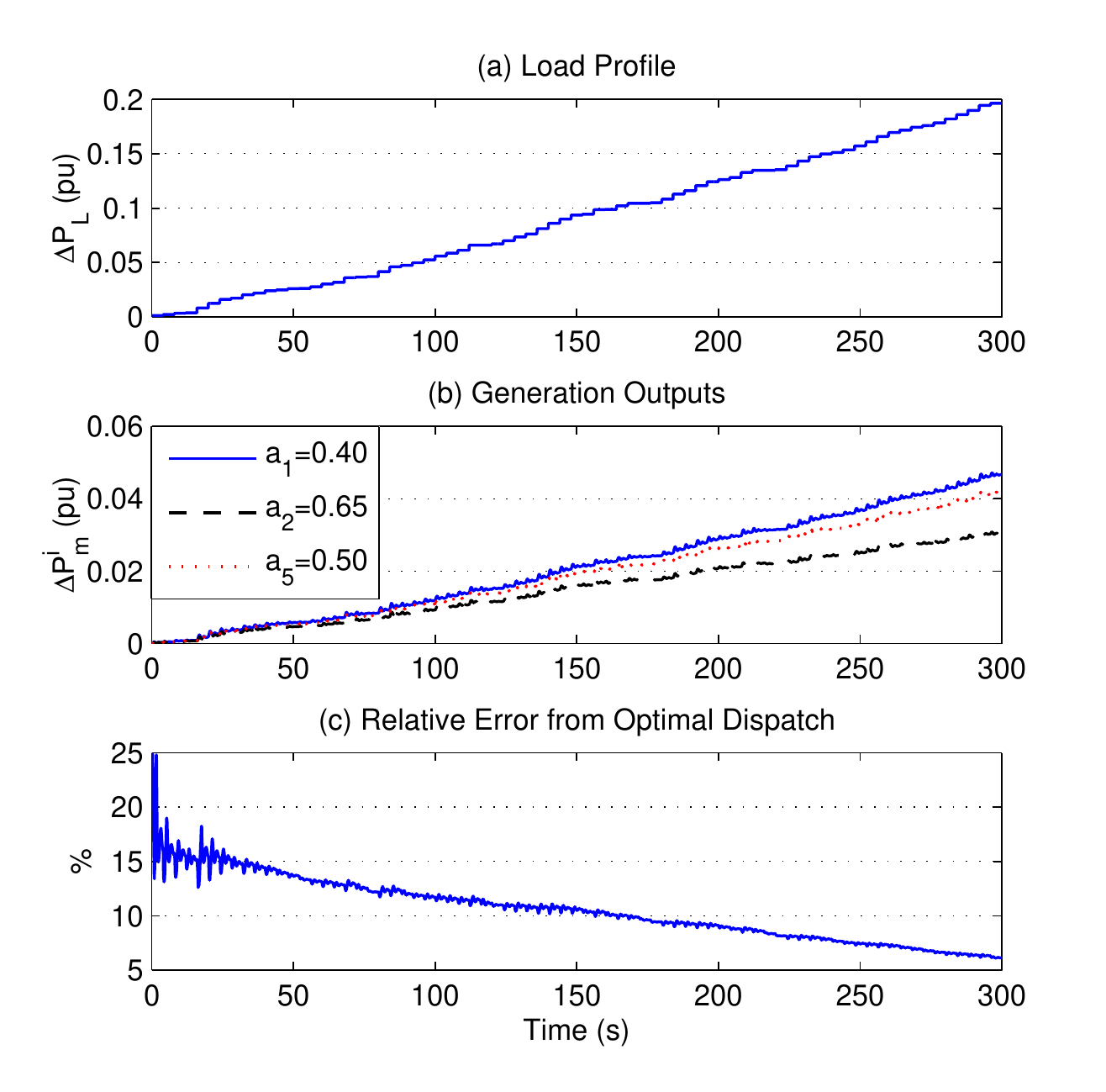}\vspace{-0.3cm}
        \caption{Cost effectiveness illustration of our distributed control.}\vspace{-0.3cm}
        \label{fig:costEffective}
        \end{center}
\end{figure}

Finally, we study the cost effectiveness of our control scheme in a 5-min time horizon in Fig. \ref{fig:costEffective}. Assuming each resource has a quadratic cost function, we know that AGC can perform the cost effective control by setting the participating factor $\alpha_i$ for resource $i$ as
\begin{equation}\label{AGC_opt}
    \alpha_i = a_i^{-1}(\textstyle\sum_{j\in\Omega}a_j^{-1})^{-1}.
\end{equation}
On the other hand, in Section \ref{analyticalEval}, we show that, by preserving the private information (i.e., the cost parameter $a_i$'s), the distributed control scheme can also asymptotically achieve the cost effectiveness. We use a monotonically increasing load, shown in Fig. \ref{fig:costEffective}(a), to better illustrate this property. The $a_i$'s for the five resources are 0.4, 0.65, 0.45, 0.6, and 0.5, respectively. The maximal relative deviation from the optimal dispatch is illustrated in Fig. \ref{fig:costEffective}(c): during the 5-min horizon, using $\Delta T=4$s, our control scheme successfully reduces the relative error from 25\% to 7\%. The large initial relative error is introduced by the global innovation term in the updating rule, and then the consensus terms dominate and reduce the relative error over time.

\vspace{-0.2cm}

\subsection{Multi-area Scenario}

\begin{figure}[t]
      \begin{center}
        \includegraphics[width=4.0cm]{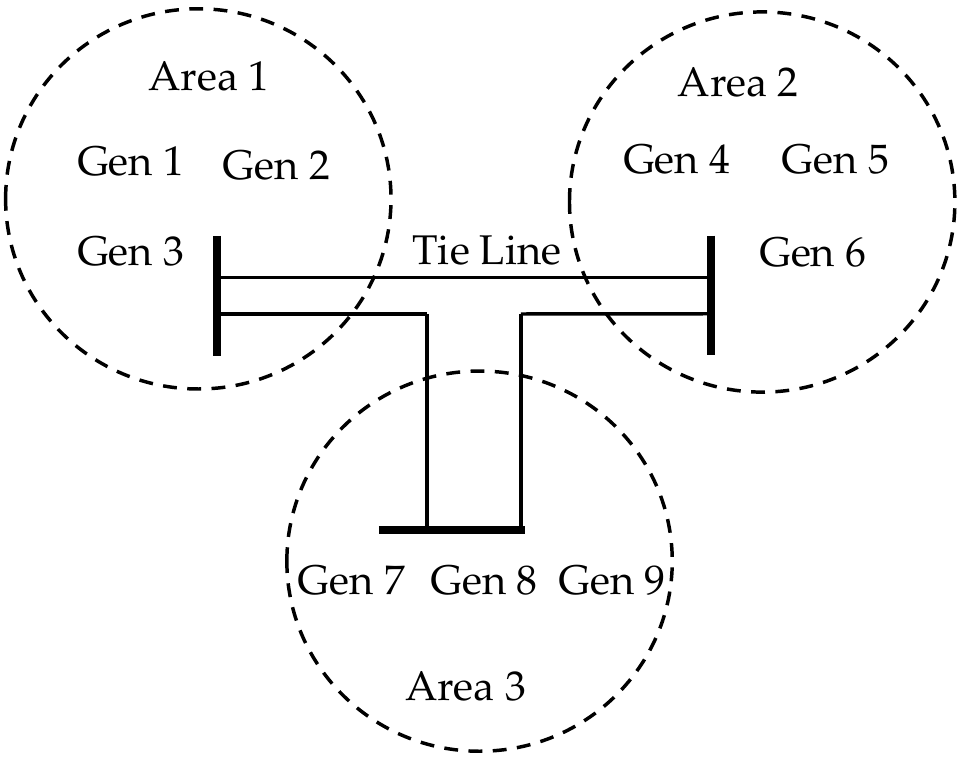}\vspace{-0.3cm}
        \caption{Three area model.}\vspace{-0.7cm}
        \label{fig:threeArea}
        \end{center}
\end{figure}

\begin{figure}[t]
      \begin{center}
        \includegraphics[width=5.6cm]{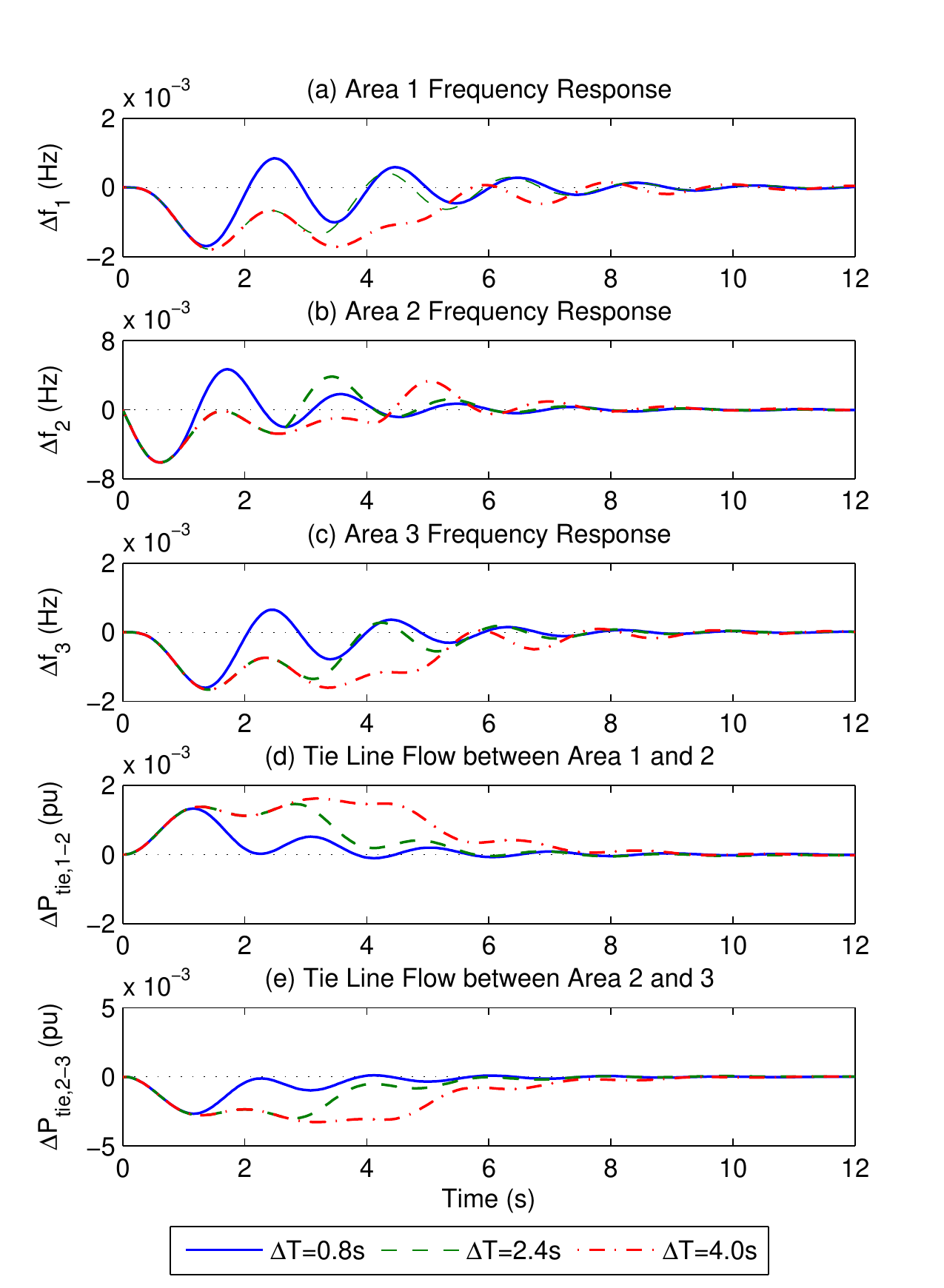}\vspace{-0.2cm}
        \caption{Frequency response with step load change in Area 2.}\vspace{-0.7cm}
        \label{fig:multiarea}
        \end{center}
\end{figure}

To verify the stability of our approach in the multi-area scenario,  we consider a three area model as shown in Fig. \ref{fig:threeArea}. Each area consists of three regulation resources. The system parameters are generated similarly as the single area scenario. 
When encountering a step load change of 0.005 pu in Area 2, the frequency responses of the three areas are illustrated in Fig. \ref{fig:multiarea}(a)-(c). The larger $\Delta T$ is, the longer our control scheme takes to stabilize the system.

The tie-line flow constraints are all set to be zero. As shown in Fig. \ref{fig:multiarea}(d)-(e), our control scheme guarantees the tie line flow constraints while stabilizing the frequency for each area.

\vspace{-0.1cm}
\section{Future Work and Conclusions} \label{conclusions}
This paper introduces a cost effective secondary frequency control framework, which leverages peer-to-peer communication to enhance the control performance. We employ a consensus-plus-global-innovation approach to design a distributed control scheme. Theoretical analysis and simulation results further illustrate the stability and cost-effectiveness of this scheme in both single-area as well as multi-area scenarios.

This work can be extended in various directions. For instance, we would like to quantify the impact of communication delay on the distributed control scheme. Also, we want to design a fully distributed control scheme, where each resource measures the frequency information locally. This is significantly challenging, since, in such a case, the frequency throughout the system is different. The regulation resources need to first perform a consensus to obtain the `global' frequency information to avoid the resources' performing the regulation against each other.

\bibliography{distributed_control}

\appendix

\subsection{Condition for Ignoring Ramping Constraints}
\label{app1}

Under Assumption \ref{zero-time}, Theorem \ref{thm:costeffective} shows that if $\|\Delta P_L(t)-\Delta P_L(t-1)\|\le \epsilon, \ \forall t\in\mathcal{T}$, there exists a $c$, such that
$\|\lambda_n^t-\lambda^t\|\le c \epsilon.$
Assumption \ref{zero-time} further guarantees:
\begin{equation}\label{app_de_u}
\begin{aligned}
 \!\!\!\! \!\!\!\! & & & \|u_i(t+1)-u_i(t)\| \\
  \!\!\!\!\!\!\!\!  & & \le & \beta \textstyle\sum_{l\in\Omega_i} \| \lambda_i(t)\!-\!\lambda_l(t) \|\! +\!n^{-1}\|\Delta P_L(t)\!-\!\Delta P_L(t\!-\!1)\|\!\!\\
  \!\!\!\!\!\!\!\!  & & \le & \beta\textstyle \sum_{l\in\Omega_i} (\| \lambda_i(t)-\lambda(t) \| +\| \lambda(t)-\lambda_l(t) \| )+ n^{-1}\epsilon\\
  \!\!\!\! \!\!\!\! & & \le & ( 2\beta c |\Omega_i| + n^{-1}) \epsilon.
    \end{aligned}
\end{equation}
Thus, as long as
\begin{equation}\label{cond1}
    ( 2\beta c |\Omega_i| + n^{-1}) \epsilon \le r_i, \forall i\in\Omega,
\end{equation}
where $r_i$ is regulation resource $i$'s ramping limit, no ramping constraints will be binding. Hence, in such conditions, our relaxation is exact.

\end{document}